\documentclass[12pt,a4paper]{amsart}
\usepackage[latin1]{inputenc}
\usepackage{amsmath}
\usepackage{amsfonts}
\usepackage{amssymb}
\usepackage{amsthm}
\usepackage{anysize}
\usepackage{enumitem}
\usepackage{amscd}
\usepackage{hyperref}
\usepackage[square, comma, numbers, sort&compress]{natbib}
\usepackage{indentfirst}
\marginsize{3.0cm}{3.0cm}{3.0cm}{3.0cm}
\setenumerate{label={\normalfont(\arabic*)}}

\newcommand{\form}[1]{{\langle #1 \rangle }}

\newtheorem{theorem}{Theorem}[section]

\newtheorem{lemma}[theorem]{Lemma}
\newtheorem{proposition}[theorem]{Proposition}

\theoremstyle{definition}

\theoremstyle{remark}
\newtheorem{remark}[theorem]{Remark}

\numberwithin{equation}{section}
\begin{document}
\title[The Artin-Spinger Theorem for quadratic forms over semi-local rings]{The Artin-Springer Theorem for quadratic forms over semi-local rings with finite residue fields}
\author{Stephen Scully}
\date{\today}
\address{Department of Mathematical and Statistical Sciences, University of Alberta, Edmonton AB T6G 2G1, Canada}
\email{stephenjscully@gmail.com}
\thanks{The author was supported by a PIMS postdoctoral fellowship held at the University of Alberta.}

\maketitle

\begin{abstract} Let $R$ be a commutative and unital semi-local ring in which 2 is invertible. In this note, we show that anisotropic quadratic spaces over $R$ remain anisotropic after base change to any odd-degree finite \'{e}tale extension of $R$. This generalization of the classical Artin-Springer theorem (concerning the situation where $R$ is a field) was previously established in the case where \emph{all} residue fields of $R$ are \emph{infinite} by I. Panin and U. Rehmann. The more general result presented here permits to extend a fundamental isotropy criterion of I. Panin and K. Pimenov for quadratic spaces over regular semi-local domains containing a field of characteristic $\neq 2$ to the case where the ring has at least one residue field which is finite.
\end{abstract}

\section{Introduction}

Let $R$ be a commutative and unital ring in which 2 is invertible. By a (free) \emph{quadratic space} over $R$, we mean a pair $(V,q)$ consisting of a free $R$-module $V$ of finite rank and a non-degenerate quadratic form $q$ on $V$ (where, by \emph{non-degenerate}, we mean that the polar bilinear form $b_q$ of $q$ (see \S 2.A) has trivial radical, i.e., that the $R$-linear map $V \rightarrow V^*$ given by $v \mapsto b_q(v,-)$ is bijective). Given an extension of (commutative, unital) rings $R \hookrightarrow S$, there is an induced base-change procedure transforming quadratic spaces over $R$ to quadratic spaces over $S$; namely, to any quadratic space $(V,q)$ over $R$, one associates the quadratic space $(V_S,q_S)$ over $S$ obtained by setting $V_S = V \otimes_R S$ and taking $q_S$ to be the unique (non-degenerate) quadratic form on $V_S$ which restricts to $q$ on $V$ (via the canonical inclusion $V \hookrightarrow V_S$).

A quadratic space $(V,q)$ over $R$ is said to be (strictly) \emph{isotropic} if there exists a \emph{unimodular} vector $v \in V$ for which $q(v) = 0$. If no such unimodular vector exists, then we say that $(V,q)$ is (strictly) \emph{anisotropic}. In other words, $(V,q)$ is anisotropic if the \emph{projective} quadric $\lbrace q = 0 \rbrace \subseteq \mathbb{P}(V)$ over $R$ admits no $R$-valued points.

A basic, yet fundamental problem for the theory of quadratic spaces over rings is that of understanding the isotropy behaviour of quadratic spaces after base change to another ring of interest. In \cite{Panin}, I. Panin proved the following important result in this direction: If $(V,q)$ is an anisotropic quadratic space over a regular semi-local domain $R$ containing a field of characteristic $0$, then it remains anisotropic after base change to the fraction field of $R$. Later, I. Panin and K. Pimenov (\cite{PaninPimenov}) showed that the statement also holds if $R$ is a regular semi-local domain of finite characteristic (different from 2), provided that \emph{all} residue fields of $R$ are \emph{infinite}. The significance of these results is that they permit to reduce a variety of problems concerning quadratic spaces over sufficiently nice semi-local domains to the case of \emph{fields}, where a rich and extensive theory is already available. Recent applications to the study of Witt groups and Milnor-Witt $K$-groups of regular equicharacteristic semi-local domains may be found in \cite{Gille} and \cite{GilleScullyZhong}.

One of the main ingredients used in the aforementioned work of Panin and Panin-Pimenov is the following generalization of the classical Artin-Springer theorem (\cite{Springer}) concerning the existence of odd-degree points on quadrics: If $(V,q)$ is an anisotropic quadratic space over a semi-local ring $R$ with \emph{infinite} residue fields (of characteristic $\neq 2$), then $(V,q)$ remains anisotropic after base change to any odd-degree finite \'{e}tale extension of $R$ (the Artin-Springer theorem being the special case where $R$ is a \emph{field}). This statement, which was established by Panin and U. Rehmann in \cite{PaninRehmann}, was used by Panin and Pimenov to prove the following more precise version of their result (cf. \cite[Proof of Thm. 1.1]{PaninPimenov}): If $R$ is a regular semi-local domain containing a field of characteristic $\neq 2$, and the Artin-Springer theorem holds for quadratic spaces over $R$, then anisotropic quadratic spaces over $R$ remain anisotropic after base change to the fraction field of $R$.

The purpose of this note is to show that the Artin-Springer theorem holds for quadratic spaces over \emph{any} semi-local ring in which 2 is invertible (see Theorem 4.1 below). In light of the above discussion, this makes the Panin-Pimenov isotropy criterion available for quadratic spaces over any regular semi-local domain of finite characteristic $\neq 2$, which was our original objective. The stated extension of the Artin-Springer theorem is achieved by means of an elementary reduction to the case already treated in \cite{PaninRehmann}. The key tool used to carry out this reduction is Proposition 3.1 below, which is perhaps of some interest in its own right. Some additional consequences of our extension results are also briefly discussed in \S 5. \\

\noindent {\bf Terminology.} In this note, all rings are assumed to be commutative and unital. By a \emph{variety}, we mean an integral scheme of finite type over a field. \\

\noindent {\bf Acknowledgements.} I would like to thank Stefan Gille for bringing the problem discussed in this text to my attention. I would also like to thank Olivier Haution and Ivan Panin for helpful discussions.

\section{Preliminaries}

In this section we recall some facts and fix some notation concerning quadratic forms over rings (in which 2 is invertible). The main reference here is \cite{Baeza}, but we do not assume much more than a basic knowledge of the theory of quadratic forms over fields. Note that we only consider quadratic forms defined on \emph{free} modules of finite rank, as opposed to quadratic forms defined on arbitrary projective modules. In particular, the Witt rings considered below are Witt rings of \emph{free} quadratic spaces. 

\subsection{Quadratic spaces over rings} Let $R$ be a ring in which $2$ is invertible and let $V$ be an $R$-module. By a \emph{quadratic form} on $V$, we mean a map $q \colon V \rightarrow R$ such that (i) $q(\lambda v) = \lambda^2 q(v)$ for all $(\lambda,v) \in R \times V$, and (ii) the symmetric form $b_q \colon V \times V \rightarrow R$ given by $b_q(v,w) = q(v+w) - q(v) - q(w)$ is $R$-bilinear. If $W$ is an $R$-submodule of $V$, then the \emph{orthogonal complement} of $W$ in $V$ (with respect to $q$) is defined as the $R$-submodule $W^\perp$ consisting of all vectors $v \in V$ such that $b_q(v,w) = 0$ for all $w \in W$. If the $R$-linear map $W \rightarrow W^*$ given by $w \mapsto b_q(w,-)$ is bijective, then we say that $W$ is \emph{non-degenerate} (with respect to $q$). In this case, $V$ decomposes as the internal direct sum of $W$ and $W^\perp$. If $V$ itself is non-degenerate, then we say that the quadratic form $q$ is \emph{non-degenerate}. In this note, we will be exclusively concerned with non-degenerate quadratic forms on \emph{free} modules of \emph{finite rank}. To this end, it is worth noting the following:

\begin{lemma} Let $(V,q)$ be a pair consisting of an $R$-module $V$ and a quadratic form $q$ on $V$, let $W$ be an $R$-submodule of $V$, and let $n$ be a positive integer. Then $W$ is free and non-degenerate of rank $n$ if and only if there exist $n$ $R$-module generators $w_1,\hdots,w_n$ of $W$ such that the matrix $\big(b_q(w_i,w_j)\big)$ is invertible over $R$.
\begin{proof} The only point worthy of remark is that the invertibility of the matrix $M := \big(b_q(w_i,w_j)\big)$ implies the freeness of $Rw_1 + \cdots + Rw_n \subseteq V$. Observe, however, that if $r_1w_1 + \cdots + r_nw_n = 0$ for some $r_1,\hdots,r_n \in R$, then the vector $(r_1,\hdots,r_n)$ lies in the kernel of $M$, whence the claim. \end{proof} \end{lemma}

By a (free) \emph{quadratic space} over $R$, we mean a pair $(V,q)$ consisting of a free $R$-module $V$ of finite rank and a non-degenerate quadratic form $q$ on $V$. The \emph{rank} of $(V,q)$ is defined as the rank of the free $R$-module $V$. By a \emph{subspace} of $(V,q)$ we mean a quadratic space $(W,p)$ such that $W$ is a free non-degenerate submodule of $V$ and $p = q|_W$ (i.e., $p$ is the restriction of $q$ to $W$).
 
Given an $R$-algebra $S$, we can define a quadratic space $(V_S,q_S)$ over $S$ by setting $V_S = V \otimes_R S$ and taking $q_S$ to be the unique (non-degenerate) quadratic form on $V_S$ which restricts to $q$ on $V$ (we suppress from the notation the dependency of this construction on the fixed $R$-algebra structure on $S$). In the other direction, if $S$ is a \emph{free} $R$-algebra of finite rank and we are given an $R$-linear map $s \colon S \rightarrow R$ such that the quadratic form $S \rightarrow R,\;x \mapsto s(x^2)$ is non-degenerate, then, given any quadratic space $(W,p)$ over $S$, the pair $(_RW, s \circ p)$ is a quadratic space over $R$ (where by $_R W$ we mean $W$ regarded as an $R$-module via the given $R$-algebra structure on $S$).
 
If $(U,p)$ is another quadratic space over $R$, then an \emph{isometry} from $(U,p)$ to $(V,q)$ is an $R$-linear isomorphism $\phi \colon U \rightarrow V$ such that $q\big(\phi(u)\big) = p(u)$ for all $u \in U$. If such an isometry exists, then we say that $(U,p)$ and $(V,q)$ are \emph{isometric}. The set $\mathrm{O}(V,q)$ of all isometries from a fixed quadratic space $(V,q)$ to itself is equipped with a natural group structure, and is called the \emph{orthogonal group} of $(V,q)$. If $w \in V$ is such that $q(w)$ is invertible in $R$, then the map $v \mapsto v - \big(q(w)^{-1}b_q(v,w)\big)w$ is an element of $\mathrm{O}(V,q)$. Self-isometries of this kind are called \emph{reflections}.

The \emph{orthogonal sum} $(U \oplus V, p \perp q)$ and \emph{tensor product} $(U \otimes_R V, p \otimes q)$ of quadratic spaces $(U,p)$ and $(V,q)$ over $R$ are defined in the obvious way. These operations descend to the level of isometry classes, thus equipping the set of all isometry classes of quadratic spaces over $R$ with the structure of a commutative unital semi-ring $\big($with unit given by the class of $(R, x \mapsto x^2)\big)$. The Grothendieck completion of this semi-ring is called the \emph{Grothendieck-Witt ring} of (free) quadratic spaces over $R$, and is denoted $\mathrm{GW}(R)$ (see \cite[\S I.2]{Baeza} for more details).

Given an invertible element $a$ in $R$, we write $\form{a}$ for the rank-1 quadratic space $(R, x \mapsto ax^2)$ over $R$. Given $n$ invertible elements $a_1,\hdots,a_n$ in $R$, we then write $\form{a_1,\hdots,a_n}$ for the $n$-fold orthgonal sum $\form{a_1} \perp \cdots \perp \form{a_n}$. A quadratic space isometric to one of this kind is said to be \emph{diagonalizable}. Note that every rank-1 quadratic space is diagonalizable. If a quadratic space $(V,q)$ over $R$ admits a free non-degenerate rank-1 $R$-submodule $L$, then we have the orthogonal decomposition $(V,q) = (L,q|_L) \perp (L^\perp,q|_{L^\perp})$. In light of the preceding remarks, this implies:

\begin{lemma} If every quadratic space over $R$ admits a free non-degenerate rank-$1$ $R$-submodule, then every quadratic space over $R$ is diagonalizable. \end{lemma}

If $(V,q)$ is a quadratic space over $R$, then a vector $v \in V$ is said to be (strictly) \emph{isotropic} if $v$ is unimodular (i.e., $v$ generates a free rank-$1$ $R$-submodule of $V$) and $q(v) = 0$. A quadratic space which admits an isotropic vector is said to be (strictly) \emph{isotropic}, and (strictly) \emph{anisotropic} otherwise. The most important example of an isotropic quadratic space over $R$ is the \emph{hyperbolic plane} $\mathbb{H} = (R \oplus R, (x,y) \mapsto xy)$. A quadratic space is said to be \emph{hyperbolic} if it is isometric to an orthogonal sum of hyperbolic planes. The set of isometry classes of hyperbolic spaces over $R$ forms an ideal in the Grothendieck-Witt ring $\mathrm{GW}(R)$. The quotient of $\mathrm{GW}(R)$ by this ideal is called the \emph{Witt ring} of (free) quadratic spaces over $R$, and is denoted by $\mathrm{W}(R)$.

If $f \colon R \rightarrow S$ is a homomorphism of rings, then the assignment $(V,q) \mapsto (V_S,q_S)$ descends to a ring homomorphism $f^* \colon \mathrm{W}(R) \rightarrow \mathrm{W}(S)$ called \emph{restriction} along $f$. If $f$ equips $S$ with the structure of a free $R$-module of finite rank, and $s \colon S \rightarrow R$ is an $R$-linear map such that the quadratic form $S \rightarrow R,\;x \mapsto s(x^2)$ is non-degenerate, then the assignment $(W, p) \mapsto (_R W, s \circ p)$ gives rise to a ring homomorphism $s_* \colon \mathrm{W}(S) \rightarrow \mathrm{W}(R)$, called (Scharlau) \emph{transfer} along $s$. The pair ($f^*,s_*$) then satisfies the obvious projection formula; in particular, the composition
\begin{equation*} \mathrm{W}(R) \xrightarrow{f^*} \mathrm{W}(S) \xrightarrow{s_*} \mathrm{W}(R) \end{equation*}
coincides with multiplication by the element $s_*(\form{1})$ (see \cite[I.2.12]{Baeza}).

\subsection{Quadratic spaces over finite direct products of rings} Let $R_1,\hdots,R_s$ be rings in which 2 is invertible, and let $R = R_1 \times \cdots \times R_s$ be their direct product. Given a quadratic space $(V,q)$ over $R$, we obtain, for each $1 \leq i \leq s$, an induced quadratic space $(V_i,q_i) := (V_{R_i},q_{R_i})$ over $R_i$ by restriction along the natural projection $R \rightarrow R_i$. Note that the $s$ quadratic spaces $(V_1,q_1),\hdots,(V_s,q_s)$ completely determine $(V,q)$. Indeed, for each $i$, the map $R_i \hookrightarrow R$ given by $x \mapsto x \cdot 1$ gives rise to an inclusion $V_i \subseteq V$ under which $q|_{V_i} = q_i$. Via these inclusions, an arbitrary vector $v \in V$ decomposes uniquely as a sum $v = v_1 + \cdots + v_s$ with $v_i \in V_i$ for each $i$. The quadratic form $q$ then acts on $v$ by the formula $q(v) = q_1(v_1) + \cdots + q_s(v_s)$. Note, furthermore, that the subsets $V_i \subseteq V$ are pairwise mutually orthogonal with respect to the quadratic form $q$. This readily implies the following statements:

\begin{lemma} Let $R$ and $(V,q)$ be as above.
\begin{enumerate}[leftmargin=*] \item If $v \in V$, then $q(v) = 0$ if and only if $q_i(v_i) = 0$ for all $1 \leq i \leq s$.
\item Vectors $v_1,\hdots,v_n \in V$ generate a free non-degenerate rank-$n$ $R$-submodule of $V$ if and only if, for each $1 \leq i \leq s$, $(v_1)_i,\hdots, (v_n)_i$ generate a free non-degenerate rank-$n$ $R_i$-submodule of $V_i$.\end{enumerate} \end{lemma}

\subsection{Quadratic spaces over semi-local rings} The basic strategy for studying quadratic forms over general rings is that of reduction to the case of \emph{fields}, where a rich and extensive theory has been developing since its initiation by E. Witt in \cite{Witt}. As expounded in the book \cite{Baeza}, this strategy can be applied successfully to the study of quadratic forms over \emph{semi-local} rings using reduction modulo the Jacobson radical and the remarks of the previous subsection. For later reference, we now point out some basic results which emerge from this philosophy. For the remainder of this subsection, we let $R$ be a semi-local ring in which 2 is invertible and set $\overline{R} = R/\mathrm{Jac}(R)$, where $\mathrm{Jac}(R)$ denotes the Jacobson radical of $R$. Then $\overline{R}$ decomposes (essentially uniquely) as a finite direct product of fields, say $\overline{R} \simeq k_1 \times \cdots \times k_s$. If $V$ is an $R$-module, then we write $\overline{V}$ for the $\overline{R}$-module obtained by restricting $V$ along the quotient map $R \mapsto \overline{R}$. Similarly, for each $1 \leq i \leq s$, we write $\overline{V}_i$ for the $k_i$-vector space obtained from $V$ by restriction along $R \rightarrow k_i$. The image of a vector $v \in V$ under the canonical projection $V \rightarrow \overline{V}$ will be denoted by $\overline{v}$. In the same way, for each $1 \leq i \leq s$, we write $\overline{v}_i$ for the image of $v$ under the projection $V \rightarrow \overline{V}_i$. Note that we have $\overline{v} = \overline{v}_1 + \cdots + \overline{v}_s$. Furthermore, $v$ is unimodular if and only if $\overline{v}_i \neq 0$ for all $1 \leq i \leq s$. If $q$ is a quadratic form on $V$, then the induced quadratic form on $\overline{V}$ will be denoted by $\overline{q}$. Similarly, for each $1 \leq i \leq s$, the induced quadratic form on the $k_i$-vector space $\overline{V}_i$ will be denoted by $\overline{q}_i$. In particular, given a quadratic space $(V,q)$ over $R$, we obtain in this way a quadratic space $(\overline{V},\overline{q})$ over $\overline{R}$ and, for each $1 \leq i \leq s$, a quadratic space $(\overline{V}_i,\overline{q}_i)$ over $k_i$.

\begin{lemma} Let $(V,q)$ be a quadratic space over $R$. Then $v_1,\hdots,v_n \in V$ generate a free non-degenerate rank-$n$ $R$-submodule of $V$ if and only if, for each $1 \leq i \leq s$, $(\overline{v_1})_i,\hdots,(\overline{v_n})_i$ generate a non-degenerate rank-$n$ $k_i$-linear subspace of $\overline{V}_i$.
\begin{proof} In view of Lemma 2.3 (2), it suffices to show that $v_1,\hdots,v_n$ generate a free non-degenerate rank-$n$ $R$-submodule of $V$ if and only if $\overline{v_1},\hdots,\overline{v_n}$ generate a free non-degenerate rank-$n$ $\overline{R}$-submodule of $\overline{V}$. This follows from Lemma 2.1, since an element of $R$ (in this case the determinant of the matrix $\big(b_q(v_i,v_j)\big)$) is invertible in $R$ if and only if its image under the projection $R \rightarrow \overline{R}$ is invertible in $\overline{R}$. \end{proof} \end{lemma}

Any quadratic space over a field of characteristic $\neq 2$ admits a non-degenerate subspace of rank-$1$ (see \cite[Prop. 7.29]{EKM}). By the preceding lemma, it follows that every quadratic space over the semi-local ring $R$ admits a free non-degenerate submodule of rank 1. Lemma 2.2 therefore implies:

\begin{lemma} Every quadratic space over $R$ is diagonalizable. \end{lemma}

A slightly more subtle observation is the following important lemma:

\begin{lemma} Let $(V,q)$ be a quadratic space over $R$, and let $v,w \in V$ be unimodular vectors such that $\overline{q}_i(\overline{w}_i) = \overline{q}_i(\overline{v}_i)$ for all $1 \leq i \leq s$. Then there exists $u \in V$ such that $q(u) = q(v)$ and $\overline{u}_i = \overline{w}_i$ for all $1 \leq i \leq s$.
\begin{remark} The condition that $\overline{q}_i(\overline{w}_i) = \overline{q}_i(\overline{v}_i)$ for all $1 \leq i \leq s$ just means that $\overline{q(w)} = \overline{q(v)}$. Similarly, the assertion that $\overline{u}_i = \overline{w}_i$ for all $1 \leq i \leq s$ simply means that $\overline{u} = \overline{w}$. We choose to state the lemma this way for later reference. \end{remark}
\begin{proof} The rank-1 case is trivial, so we can assume that $(V,q)$ has rank at least 2. Since $v$ and $w$ are unimodular, we have $\overline{v}_i \neq 0 \neq \overline{w}_i$ for all $1 \leq i \leq s$. Thus, by an observation of E. Witt (see \cite[\S II.8]{EKM}), there exists, for each $1 \leq i \leq s$, an isometry $\tau_i \in \mathrm{O}(\overline{V}_i,\overline{q}_i)$ which is a product of (at most 2) reflections and which satisfies $\tau_i(\overline{v}_i) = \overline{w}_i$. Since $\mathrm{rank}(V) \geq 2$, every reflection in $\mathrm{O}(\overline{V}_i,\overline{q}_i)$ extends to a reflection in $\mathrm{O}(\overline{V},\overline{q})$ which fixes $\overline{v}_j$ for all $j \neq i$. It follows that there exists $\tau \in \mathrm{O}(\overline{V},\overline{q})$ which is a product of reflections and which sends $\overline{v}$ to $\overline{w}$. Since an element of $R$ is invertible if and only if its image in $\overline{R}$ is invertible, every reflection in $\mathrm{O}(\overline{V},\overline{q})$ lifts to a reflection in $\mathrm{O}(V,q)$. As a result, $\tau$ lifts to an isometry $\phi \in \mathrm{O}(V,q)$. The vector $u := \phi(v)$ then has the desired properties. \end{proof} \end{lemma}

\begin{remark} Using the Cartan-Dieudonn\'{e} theorem on the generation of the orthogonal group of a quadratic space over a \emph{field} by reflections (see \cite[Theorem I.5.4]{Scharlau}), one can show by similar means that the natural projection $\mathrm{SO}(V,q) \rightarrow \mathrm{SO}(\overline{V},\overline{q})$ on \emph{special} orthogonal groups is surjective. This observation is due to M. Knebusch (see \cite{Knebusch}). If $R$ is a \emph{local} ring, then even $\mathrm{O}(V,q) \rightarrow \mathrm{O}(\overline{V},\overline{q})$ is surjective. \end{remark}

It is worth stating explicitly the following special case of Lemma 2.6 which will be needed later:

\begin{lemma} Let $(V,q)$ be an isotropic quadratic space over $R$. Suppose that, for each $1 \leq i \leq s$, we are given a non-zero vector $\mathbf{v}_i \in \overline{V}_i$ such that $\overline{q}_i(\mathbf{v}_i) = 0$. Then there exists a (strictly) isotropic vector $v \in V$ such that $\overline{v}_i = \mathbf{v}_i$ for all $1 \leq i \leq s$. \end{lemma}

Now, using Lemma 2.6 (see also Remark 2.8), one can readily show that Witt's cancellation and decomposition theorems hold for (free) quadratic spaces over $R$ (see \cite[Ch. III]{Baeza}). As a result, the standard characterizations of isotropy for quadratic spaces over fields are all valid over the semi-local ring $R$:

\begin{lemma} Let $(V,q)$ be a quadratic space over $R$ and let $Q \subseteq \mathbb{P}(V)$ be the projective $R$-quadric defined by the vanishing of $q$. Then the following are equivalent:
\begin{enumerate}[leftmargin=*] \item $(V,q)$ is isotropic.
\item $(V,q) \simeq \mathbb{H} \perp (W,p)$ for some quadratic space $(W,p)$ over $R$.
\item There exists a quadratic space $(W,p)$ over $R$ with $\mathrm{rank}(W) < \mathrm{rank}(V)$ such that $(V,q) = (W,p)$ in $\mathrm{W}(R)$.
\item $Q$ has an $R$-valued point.
\item $Q$ is rational over $R$.
\item The affine quadric $\lbrace q = 0 \rbrace \subseteq \mathbb{A}(V)$ is rational over $R$. \end{enumerate}
\begin{proof} If (2) holds, then we can choose coordinates $(x:y:z_1:\cdots:z_n)$ on $V^*$ so that $q(x,y,z_1,\hdots,z_n) = xy + p(z_1,\hdots,z_n)$. The open complement in $Q$ of the principal divisor $\lbrace x = 0 \rbrace$ is then isomorphic to $\mathbb{A}_R^n$, and so (5) holds. The remaining implications are either trivial or straightforward consequences of Lemma 2.6 and standard facts concerning quadratic forms over fields; for the sake of brevity we omit the details and refer the reader instead to \cite[Ch. III]{Baeza}.\end{proof} \end{lemma}

\section{Key proposition}

The key ingredient needed to extend the Artin-Springer theorem to quadratic spaces over arbitrary semi-local rings in which 2 is invertible is the ($n=3$ case of the) following proposition, which is of some interest in its own right:

\begin{proposition} Let $R$ be a semi-local ring in which 2 is invertible, and let $S$ be a degree-$n$ finite \'{e}tale extension of $R$ for some (not necessarily odd) integer $n \geq 2$. If $(V,q)$ is a quadratic space of rank $>n$ over $R$ such that $(V_S,q_S)$ is isotropic, then $(V,q)$ contains a rank-$n$ subspace which becomes isotropic over $S$. \end{proposition}

The case where $n=2$ is well known and can be made rather more precise (see \cite[Thm. V.4.2]{Baeza}). The general case of Proposition 3.1 follows immediately from the more specific Lemma 3.2 below. Before proceeding, let us pause to fix some notation and terminology which will also be used in the next section:\\

\noindent {\bf Notation and Terminology.} Let $R$ and $S$ be as in the statement of Proposition 3.1. By \cite[18.4.5]{EGAIV}, we can write $S = R[t]/\big(f(t)\big)$, where $f$ is a monic separable polynomial of degree $n$ over $R$. By \emph{monic}, we mean that the leading coefficient of $f$ is invertible in $R$; by \emph{separable}, we mean that the discriminant of $f$ is invertible in $R$. Both conditions are local in the following sense: a polynomial over $R$ is monic (resp. separable) if and only if it is monic (resp. separable) when restricted to every residue field of $R$. The image of the variable $t$ in $S$ will be denoted by $\theta$. Note that $S$ is a free rank-$n$ $R$-module with basis $1, \theta,\hdots,\theta^{n-1}$.\\

Our refinement of Proposition 3.1 may now be stated as follows:

\begin{lemma} In the situation of Proposition 3.1, there exist vectors $v_0,v_1,\hdots,v_{n-1} \in V$ such that
\begin{enumerate}[leftmargin=*] \item $v_0 + \theta v_1 + \cdots + \theta^{n-1}v_{n-1} \in V_S$ is a $($strictly$)$ isotropic vector for $q_S$. 
\item $v_0,v_1,\hdots,v_{n-1}$ generate a free non-degenerate rank-$n$ submodule of $V$. \end{enumerate}

\begin{proof} Let us first remark that the unimodularity of the vector in (1) is guaranteed by condition (2), and so we will suppress this point in the discussion that follows.

We begin the proof with some reductions. Firstly, since $S$ is a finite \'{e}tale extension of $R$, its Jacobson radical is generated (as an ideal of $S$) by the Jacobson radical of $R$. Thus, by Lemmas 2.4 and 2.9, the problem reduces immediately to the case where $R$ is a field, which we will now choose to denote by $k$. Before proceeding, we make one further reduction to the case where the $u$-invariant of $k$ is $\leq 2$ (that is, every quadratic space of rank $\geq 3$ over $k$ is isotropic). If $k$ is finite, then it has $u$-invariant $2$, and so no reduction is necessary (see \cite[Ex. 36.2 (3)]{EKM}). To treat the case where $k$ is infinite, consider the \emph{affine} quadric $X_{q_S} = \lbrace q_S = 0 \rbrace \subseteq \mathbb{A}(V_S)$. Since $q_S$ admits a (strictly) isotropic vector, $X_{q_S}$ is a rational $S$-scheme (Lemma 2.10). The Weil restriction $Y := R_{S/k}(X_{q_S})$ is therefore a rational $k$-variety. Note that $Y$ is a closed subvariety of $\mathbb{A}(V^{\oplus n})$: if $A$ is a $k$-algebra, then an $n$-tuple of vectors $(v_0,v_1,\hdots,v_{n-1}) \in (V_A)^{\oplus n}$ lies in $Y(A)$ if and only if $q_{S \otimes_k A}(\sum_{i=0}^{n-1}\theta^i v_i) = 0$. Consider now the $k$-morphism $g \colon \mathbb{A}(V^{\oplus n}) \rightarrow \mathbb{A}^1_k$ given by $(v_0,v_1,\hdots,v_{n-1}) \mapsto \mathrm{det}\big(b_q(v_i,v_j)\big)$. By Lemma 1.1, $n$-vectors $v_0,v_1,\hdots,v_{n-1}$ in $V$ generate a non-degenerate rank-$n$ subspace of $V$ if and only if $g(v_0,v_1,\hdots,v_{n-1}) \neq 0$. Thus, in light of the above remarks, the lemma holds if and only if $U(k) \neq \emptyset$, where $U$ is the open subvariety of $Y$ given by $U = Y \cap g^{-1}(\mathbb{A}^1_k \setminus \lbrace 0 \rbrace)$. Now $Y$ is a rational $k$-variety, and so if $k$ is infinite, then $U(k) \neq \emptyset$ if and only if $U$ is non-empty \emph{as a scheme}. In particular, in order to prove the lemma in this case, we may assume that $k$ is algebraically closed. Since any quadratic space of rank 2 or more over an algebraically closed field is evidently isotropic, we have made the necessary reduction.

Working under the assumption that the $u$-invariant of $k$ is $\leq 2$, we will now show that there exist vectors $v_0,v_1,\hdots,v_{n-1} \in V$ satisfying conditions (1) and (2) in the statement of the lemma. Let $s = [\frac{n}{2}]$. Note that $(V,q)$ has rank strictly greater than $2s$ by assumption. Since every rank-$3$ quadratic space over $k$ is isotropic, it follows that $V$ admits a non-degenerate subspace $W$ of dimension $2s + 1$ such that
$$ (W,q|_W) \simeq \underbrace{\mathbb{H} \perp \cdots \perp \mathbb{H}}_{s \text{ times}} \perp \form{a} $$
for some $a \in k^*$ (see Lemma 2.10). More explicitly, this means that we can find vectors $e_1,f_1,e_2,f_2,\hdots,e_s,f_s$ and $w$ in $V$ such that
\begin{itemize} \item $q(e_i) = q(f_i) = 0$ for all $i$,
\item $q(w) = a$,
\item $b_q(e_i,e_j) = b_q(f_i,f_j) = 0$ for all $i \neq j$,
\item $b_q(e_i,f_j) = \delta_{ij}$ for all $i,j$ (where $\delta$ is the Kronecker delta), and
\item $b_q(e_i,w) = b_q(f_i,w) = 0$ for all $i$. \end{itemize}
(the $(e_i,f_i)$ are \emph{hyperbolic pairs} in the sense of \cite[p. 40]{EKM}) To finish the proof, let us now separate the case where $n$ is even from that where $n$ is odd:\\

\noindent {\it Case 1.} $n=2s$. Suppose first that $s =1$, and let $\alpha, \beta \in k$ be such that $\theta^2 = \alpha + \beta \theta$ in $S$. A direct calculation shows that if we set
\begin{equation*} v_0 = (-\alpha - \frac{\beta^2}{4})ae_1 + f_1 -\frac{\beta}{2}w \hspace{1cm} \text{and} \hspace{1cm} v_1 = w, \end{equation*}
then $q_S(v_0 + \theta v_1) = 0$, i.e., condition (1) of the statement holds for this choice of $v_0$ and $v_1$. At the same time, the determinant of the Gram matrix $\big(b_q(v_i,v_j)\big)$ is in this case equal to $-a^2(\beta^2 + 4\alpha)$, which is non-zero because $S/k$ is \'{e}tale. Thus, condition (2) also holds for this choice of $v_0$ and $v_1$.

Suppose now that $s >1$. Let $\lambda_1,\hdots,\lambda_s \in k^*$ be such that $\sum_{i=1}^s \lambda_i = 0$, and set
$$ v_i = \begin{cases} \lambda_{i+1} e_{i+1} & \text{if } 0 \leq i \leq s-1 \\
                       f_{n-i} & \text{if } s \leq i \leq n-1 \end{cases} $$
Then, by construction, we have
$$ q_S\bigg(\sum_{i=0}^{n-1}\theta^i v_i \bigg) = \sum_{i=0}^{n-1}q(v_i)\theta^{2i} + \sum_{i<j} b_q(v_i,v_j) \theta^{i+j} = \sum_{i=1}^s \lambda_i \theta^{n-1} = 0, $$
and so the $v_i$ satisfy the first condition of the statement. Condition (2) is also satisfied, since the restriction of $q$ to the subspace spanned by the $v_i$ is a hyperbolic form. \\

\noindent {\it Case 2.} $n=2s + 1$. Let $\lambda_1,\hdots,\lambda_s \in k^*$ be such that $\sum_{i=1}^s \lambda_i = -a$, and set
$$ v_i = \begin{cases} \lambda_{i+1} e_{i+1} & \text{if } 0 \leq i \leq s-1 \\
                       w & \text{if } i = s \\
                       f_{n-i} & \text{if } s+1 \leq i \leq n-1 \end{cases} $$
Then, by construction, we have
$$ q_S\bigg(\sum_{i=0}^{n-1}\theta^i v_i \bigg) = \sum_{i=0}^{n-1}q(v_i)\theta^{2i} + \sum_{i<j} b_q(v_i,v_j) \theta^{i+j} = q(w)\theta^{n-1} + \sum_{i=1}^s \lambda_i \theta^{n-1} = (a-a)\theta^{n-1} = 0, $$
and so the $v_i$ satisfy condition (1) of the statement. Again, condition (2) is also satisfied, since the subspace of $V$ spanned by the $v_i$ is nothing else but $W$. \end{proof} \end{lemma}

\section{The Artin-Springer theorem} 

We are now ready to prove our main result:

\begin{theorem} Let $R$ be a semi-local ring in which 2 is invertible. If $(V,q)$ is an anisotropic quadratic space over $R$, then $(V_S,q_S)$ is anisotropic for any odd-degree finite \'{e}tale extension $S$ of $R$. \end{theorem}

As discussed in \S 1, the particular case where $R$ is a field is due (independently) to E. Artin (unpublished) and T.A. Springer (\cite{Springer}), while the case where \emph{all} residue fields of $R$ are \emph{infinite} was proved in \cite{PaninRehmann} by I. Panin and U. Rehmann. Our purpose is to treat the case where at least one residue field of $R$ is finite. We will argue by reduction to the case where all residue fields are ``large enough''. The key ingredient here is (the $n=3$ case of) Proposition 3.1 above. We continue with the notation of \S 3, but, as per the statement of Theorem 4.1, we now assume that $n$ is \emph{odd}. We present the proof of Theorem 4.1 in a series of short steps:

\subsection{The hyperbolicity analogue of Theorem 4.1} Let $i$ denote the natural inclusion of $R$ into $S$. Below, we will make use of the following hyperbolicity analogue of Theorem 4.1, which is well known and straightforward to prove (see \cite[\S 3]{OjangurenPanin}):

\begin{proposition} The natural restriction homomorphism $i^* \colon \mathrm{W}(R) \rightarrow \mathrm{W}(S)$ is injective. In other words, if $(V,q)$ is a quadratic space over $R$ such that $(V_S,q_S)$ is hyperbolic, then $(V,q)$ is hyperbolic.
\begin{proof} Let $s \colon S \rightarrow R$ be the $R$-linear map defined by setting $s(\theta^{n-1}) = 1$ and $s(\theta^i) = 0$ for all $0 \leq i \leq n-2$. A straightforward calculation shows that the quadratic form $S \rightarrow R, x \mapsto s(x^2)$ is isometric to $\form{1} \perp \frac{n-1}{2} \cdot \mathbb{H}$ (here $\frac{n-1}{2} \cdot \mathbb{H}$ denotes the orthogonal sum of $\frac{n-1}{2}$ copies of $\mathbb{H}$). In particular, it is non-degenerate, and so we can consider the Scharlau transfer $s_* \colon \mathrm{W}(S) \rightarrow \mathrm{W}(R)$. Recall from \S 2.1 above that the composition $\mathrm{W}(R) \xrightarrow{i^*} \mathrm{W}(S) \xrightarrow{s_*} \mathrm{W}(R)$ coincides with multiplication by the element $s_*(\form{1})$. However, by the above remarks, we have $s_*(\form{1}) = \form{1}$ in $\mathrm{W}(R)$, and so we see that $s_*$ is a section of $i^*$. Hence $i^*$ is injective, as claimed. \end{proof} \end{proposition}

\subsection{Proof of Theorem 4.1 for quadratic spaces of rank $\leq 3$} Using Proposition 4.2, we now show that Theorem 4.1 holds for quadratic spaces of rank $\leq 3$.

\begin{proposition} Theorem 4.1 holds when $(V,q)$ has rank $\leq 3$.
\begin{proof} For quadratic spaces of rank 1, the statement is trivial. Next, an isotropic space of rank 2 over a semi-local ring is hyperbolic (see Lemma 2.10), and so the rank-2 case of the statement follows from Proposition 4.2. Assume now that $(V,q)$ has rank 3. For ease of notation, we simply denote the pair $(V,q)$ by $q$. Now, by Lemma 2.5, we can write $q \simeq \form{a,b,c}$ for some $a,b,c \in R^*$. Then the form $\pi = q \perp \form{abc} = \form{a,b,c,abc}$ is similar to (i.e., isometric to a unit multiple of) the 2-fold Pfister form $\form{1,ab} \otimes \form{1,ac}$. By the main property of Pfister forms (see \cite[Cor. IV.3.2]{Baeza}), $\pi_S$ is either anisotropic or hyperbolic. Suppose that $\pi_S$ is hyperbolic. Then, by Proposition 4.2, $\pi$ is already hyperbolic, and so $q = -\form{abc}$ in the Witt ring $\mathrm{W}(R)$. But, by Lemma 2.10, this implies that $q$ is isotropic, thus contradicting our initial hypothesis. We can therefore conclude that $\pi_S$ is anisotropic. Since $q_S$ is a restriction of $\pi_S$, it must be anisotropic as well. \end{proof} \end{proposition}

\subsection{Proof of Theorem 4.1 for extensions of degree 3} Next, we use Proposition 4.3 to show that Theorem 4.1 holds in the case where $S$ has degree $3$ (with no restriction on the rank of $(V,q)$):

\begin{proposition} Theorem 4.1 holds when $S$ is a degree-3 extension of $R$.
\begin{proof} The case where $(V,q)$ has rank $\leq 3$ was treated in Proposition 4.3 above. If $(V,q)$ has rank strictly greater than $3$, then Proposition 4.3 at least implies that every rank-3 subspace of $(V,q)$ remains anisotropic over $S$. In view of the $n=3$ case of Proposition 3.1, however, this is all that is needed to make the desired conclusion. \end{proof} \end{proposition}

\subsection{Proof of Theorem 3.1 in the general case} We now conclude the proof of Theorem 3.1. The main ingredient, together with Proposition 4.4, is the following statement due (essentially) to I. Panin and U. Rehmann:

\begin{lemma}[Panin-Rehmann] There exists a natural number $M_n$ for which the following holds: Assume that all residue fields of $R$ contain at least $M_n$ elements. Then, for any quadratic space $(V,q)$ over $R$ such that $(V_S,q_S)$ is isotropic, there exist $v_0,v_1,\hdots,v_{n-1} \in V$ with the following properties:
\begin{enumerate}[leftmargin=*] \item $v_0 + \theta v_1 + \cdots + \theta^{n-1}v_{n-1} \in V_S$ is a (strictly) isotropic vector for $q_S$.
\item $q_{R[t]}(v_0 + tv_1 + \cdots + t^{n-1}v_{n-1}) \in R[t]$ is monic and separable of degree $2n-2$ $($here $R[t]$ denotes the polynomial ring in a single variable $t$ over $R)$. \end{enumerate} 
\begin{proof} This is a mild extension of \cite[Prop. 1.1]{PaninRehmann}. In fact, the statement is proved in \cite{PaninRehmann} under the stronger hypothesis that $R$ is \emph{local} with \emph{infinite} residue field. Note, however, that since $S$ is a finite \'{e}tale extension of $R$, its Jacobson radical is generated (as an ideal of $S$) by that of $R$. Thus, in view of Lemma 2.9 and the fact that condition (2) can be checked locally (see the remarks at the beginning of this section), everything reduces immediately to the case where $R$ is a field (regardless of whether $R$ is local or not). Moreover, under the assumption that $R$ is a field, it was shown in \cite[\S 3]{PaninRehmann} that there exists a non-empty rational $R$-variety $U$ with property that $U(R) \neq \emptyset$ if and only if the conclusion of the lemma holds. This conclusion is therefore valid provide that $R$ contains ''enough'' elements; the term ``enough'' is made formal here by the integer $M_n$ in the statement of the lemma (in particular, it is not necessary for $R$ to be infinite). \end{proof} \end{lemma}

We now give the proof of Theorem 4.1:

\begin{proof}[Proof of Theorem 4.1] Suppose first that at least one residue field of $R$ is finite. Let $h(x) \in R[x]$ be a monic separable polynomial of degree 3 in a single variable $x$ over $R$ which is irreducible when restricted to any finite residue field of $R$. Then $\widetilde{R} := R[x]/\big(h(x)\big)$ is a degree-3 \'{e}tale extension of $R$ such that the minimal cardinality of the finite residue fields of $\widetilde{R}$ is strictly greater than that of the finite residue fields of $R$. Now, in order to show that $(V_S,q_S)$ is anisotropic, it is sufficient to show that $(V_{\widetilde{S}},q_{\widetilde{S}})$ is anisotropic, where $\widetilde{S} = S \otimes_R \widetilde{R}$. By Proposition 4.4, however, $(V_{\widetilde{R}},q_{\widetilde{R}})$ is anisotropic. We can therefore replace $R$ by $\widetilde{R}$, and, by repeating this process sufficiently many times, we ultimately reduce to the case where all residue fields of $R$ contain at least $M_n$ elements (where $M_n$ is the natural number in the statement of Lemma 4.5). Having made this reduction, the proof proceeds as in \cite{PaninRehmann}. For the reader's convenience, we give the details. Thus, suppose for the sake of contradiction that $(V_S,q_S)$ is isotropic, and let $v_0,v_1,\hdots,v_n$ be as in Lemma 4.5 so that
\begin{enumerate}\item $v_0 + \theta v_1 + \cdots + \theta^{n-1}v_{n-1} \in V_S$ is a (strictly) isotropic vector for $q_S$.
\item $q_{R[t]}(v_0 + tv_1 + \hdots t^{n-1}v_{n-1}) \in R[t]$ is monic and separable of degree $2n-2$. \end{enumerate}
Let $v(t) = v_0 + tv_1 + \hdots t^{n-1}v_{n-1} \in V \otimes_R R[t]$. By (1), the polynomial $q_{R[t]}\big(v(t)\big) \in R[t]$ is divisible by $f(t)$, say
\begin{equation*} q_{R[t]}\big(v(t)\big) = f(t)g(t) \end{equation*}
with $g(t) \in R[t]$. By (2), $g(t)$ is monic and separable of degree $n-2$ over $R$. Let $S' = S[t]/\big(g(t)\big)$. Then $S'$ is a degree-$(n-2)$ \'{e}tale extension of $R$. We claim that $(V_{S'},q_{S'})$ is isotropic. Note first that if $\theta'$ denotes the image of $t$ in $S'$, then we have
\begin{equation*} q_{S'}\big(v(\theta')\big) = f(\theta')g(\theta') = 0 \end{equation*}
in $S'$. To prove the claim, it therefore suffices to check that $v(\theta') \in V_{S'}$ is unimodular. We claim that this follows from the separability part of (2). Since unimodularity of vectors and separability of polynomials can be checked locally (see \S 2.3 and the remarks at the beginning of \S 3), we are, to this end, free to assume that $R$ is a field. Then $g(t)$ decomposes into a product of pairwise coprime irreducible polynomials over $R$, say $g(t) = g_1(t)\cdots g_k(t)$. If $v(\theta')$ were not unimodular, then $v(t)$ would be divisible by some $g_i(t)$, whence $q_{R[t]}$ would be divisible by $g_i(t)^2$ in $R[t]$. As this would contradict the separability part of (2), the claim follows. Thus, $(V_{S'},q_{S'})$ is isotropic. Now $S'$ is an odd-degree \'{e}tale extension of $R$ which has strictly smaller degree than $S$. By repeating the above argument sufficiently many times, we come to the conclusion that $(V,q)$ is isotropic, thus contradicting our original hypothesis. This completes the argument. \end{proof}

\section{The Panin-Pimenov isotropy criterion and applications}

\subsection{The Panin-Pimenov isotropy criterion} As an immediate application of our extension of the Artin-Springer theorem (Theorem 4.1 above), we obtain the following extension of an important result due to I. Panin and K. Pimenov:

\begin{theorem} Let $R$ be a regular semi-local domain containing a field of characteristic $\neq 2$. If $(V,q)$ is an anisotropic quadratic space over $R$, then $(V_K,q_K)$ is anisotropic, where $K$ denotes the fraction field of $R$.  \end{theorem}

This theorem was originally proved by I. Panin in \cite{Panin} under the stronger assumption that $R$ contains a field of characteristic $0$. In \cite{PaninPimenov}, I. Panin and K. Pimenov treated the case where $R$ has finite characteristic different from 2, but under the assumption that \emph{all} residue fields of $R$ are \emph{infinite}. In order to extend the statement to the case where some residue field of $R$ is finite (which was our original goal), the arguments presented in \cite{PaninPimenov} showed that it would suffice to prove that the Artin-Springer theorem holds for quadratic spaces over $R$. This extension therefore follows from our Theorem 4.1 above. In more detail:

\begin{proof}[Proof of Theorem 4.1] Let $(V,q)$ be a quadratic space over $R$. If $(V_K,q_K)$ is isotropic, then it is shown in \cite[Proof of Thm. 1.1]{PaninPimenov} that $(V,q)$ becomes isotropic over an odd-degree finite \'{e}tale extension of $R$ (in fact, over $R \otimes_k k'$, where $k$ denotes the prime subfield of $R$ and $k'$ is some finite extension of $k$ of degree prime to $2$ and the characteristic of $k$). By Theorem 4.1, this implies that $(V,q)$ is isotropic. \end{proof}

\subsection{Applications} The Panin-Pimenov isotropy criterion is of fundamental importance for the study of quadratic forms over (nice) semi-local rings, since it permits to reduce many problems to the case of \emph{fields}, where an extensive literature is already available (see, e.g., \cite{EKM}). Let $R$ be a regular semi-local domain containing a field of characteristic $\neq 2$. To conclude this note, we now mention (without proof) some consequences of our extension of the Panin-Pimenov criterion to the case where at least one residue field of $R$ is finite (Theorem 5.1). These statements were known previously in the case where the residue fields of $R$ are all infinite.

\begin{itemize}[leftmargin=*] \item All results proved in \cite[\S 3]{Gille} are valid for quadratic spaces over $R$; these include the \emph{linkage theorem} (analogous to the original result of R. Elman and T. Y. Lam concerning quadratic spaces over \emph{fields}) as well as certain other purity-type statements for Pfister forms.
\item If, moreover, $R$ \emph{local}, then the purity theorem of \cite[Cor. 1]{PaninPimenov} holds for quadratic spaces over $R$, i.e., unramified quadratic spaces over $K := \mathrm{Frac}(R)$ are defined over $R$ (on the level of quadratic spaces themselves, and not only Witt groups).
\item If $R$ is \emph{local}, then the main results of \cite{ChernousovPanin1} and \cite{ChernousovPanin2} (concerning purity statements for $G_2$-torsors and certain $F_4$-torsors, respectively) are valid over $R$.
\end{itemize}

\bibliographystyle{alphaurl}

\begin{thebibliography}{EGAIV}

\bibitem[Ba78]{Baeza}
R.~Baeza.
\newblock {\em Quadratic forms over semilocal rings}.
\newblock Lecture Notes in Mathematics, Vol. 655, Springer-Verlag, Berlin-New York, 1978, vi+199 pp.

\bibitem[CP07]{ChernousovPanin1}
V.~Chernousov and I. Panin.
\newblock {\em Purity of $G_2$-torsors}.
\newblock C. R. Math. Acad. Sci. Paris {\bf 345} (2007), no.6, 307--312.

\bibitem[CP13]{ChernousovPanin2}
V.~Chernousov and I. Panin.
\newblock {\em Purity for Pfister forms and $F_4$-torsors with trivial $g_3$ invariant}.
\newblock J. Reine Angew. Math. {\bf 685} (2013), 99--104.

\bibitem[EGAIV]{EGAIV}
A.~Grothendieck.
\newblock{\em El\'{e}ments de g\'{e}ometrie alg\'{e}brique. IV.}.
\newblock Inst. Hautes \'{E}tudes Sci. Publ. Math.

\bibitem[EKM08]{EKM}
R.~Elman, N.~Karpenko, and A.~Merkurjev.
\newblock {\em The algebraic and geometric theory of quadratic forms}.
\newblock AMS Colloquium Publications, {\bf 56}, American Mathematical Society, 2008.

\bibitem[G15]{Gille}
S.~Gille.
\newblock {\em On quadratic forms over semilocal rings}.
\newblock Preprint, 2015.

\bibitem[GSZ16]{GilleScullyZhong}
S.~Gille, S.~Scully and C.~Zhong.
\newblock{\em Milnor-Witt $K$-groups of local rings}.
\newblock Adv. Math. {\bf 286} (2016), 729--753.

\bibitem[Kn69]{Knebusch}
M.~Knebusch.
\newblock {\em Isometrien \"{u}ber semilokalen Ringen. (German)}.
\newblock Math. Z. {\bf 108}, 1969, 255--268.

\bibitem[OP99]{OjangurenPanin}
M.~Ojanguren and I.~Panin.
\newblock {\em A purity theorem for the Witt group}.
\newblock Ann. Sci. \'{E}cole Norm. Sup. (4) {\bf 32} (1999), no. 1, 71--86.

\bibitem[P09]{Panin}
I.~Panin.
\newblock {\em Rationally isotropic quadratic spaces are locally isotropic}.
\newblock Invent. Math. {\bf 176} (2009), no. 2, 397--403.

\bibitem[PP10]{PaninPimenov}
I.~Panin and K.~Pimenov.
\newblock {\em Rationally isotropic quadratic spaces are locally isotropic: II}.
\newblock Doc. Math. 2010, Extra volume: Andrei A. Suslin sixtieth birthday, 515--523.

\bibitem[PR07]{PaninRehmann}
I.~Panin and U.~Rehmann.
\newblock {\em A variant of a theorem by Springer}.
\newblock Algebra i Analiz {\bf 19} (2007), no. 6, 117--125; translation in St. Petersburg Math. J. 19 (2008), no. 6, 953--959.

\bibitem[Sc85]{Scharlau}
W.~Scharlau.
\newblock {\em Quadratic and hermitian forms}.
\newblock Springer-Verlag, Berlin, 1985.

\bibitem[Sp52]{Springer}
T.A.~Springer.
\newblock {\em Sur les formes quadratiques d'indice z\'{e}ro}.
\newblock C. R. Acad. Sci. Paris 234, (1952). 1517-1519.

\bibitem[Wi37]{Witt}
E.~Witt.
\newblock {\em Theorie der quadratischen Formen in beliebigen K\"orpern}.
\newblock J. Reine Angew. Math. {\bf 176} (1937), 31--44.

\end{thebibliography}

\end{document}